\documentclass{amsart}

\usepackage{amssymb}
\usepackage{graphicx}
\usepackage[english]{babel}
\usepackage[latin1]{inputenc}

\input xy
\xyoption{all}

\theoremstyle{plain}
\newtheorem{my_thm}{Theorem}
\newtheorem{my_thm2}{Theorem}
\newtheorem{my_prop}[my_thm]{Proposition}
\newtheorem{my_prop2}[my_thm2]{Proposition}
\newtheorem{my_def}{Definition}[section]
\newtheorem{my_cor}[my_def]{Corollary}
\newtheorem{my_lem}[my_def]{Lemma}
\newtheorem{theorem}[my_def]{Theorem}
\newtheorem{proposition}[my_def]{Proposition}

\theoremstyle{remark}
\newtheorem{my_rem}[my_def]{Remark}
\newtheorem{my_ex}[my_def]{Example}

\def\Hom{\mathop{\rm Hom}\nolimits}

\def\End{\mathop{\rm End}\nolimits}
\def\Ext{\mathop{\rm Ext}\nolimits}
\def\soc{\mathop{\rm soc}\nolimits}
\def\pd{\mathop{\rm pd}\nolimits}
\def\id{\mathop{\rm id}\nolimits}

\def\mod{\mathop{\rm mod}\nolimits}
\def\top{\mathop{\rm top}\nolimits}
\def\ind{\mathop{\rm ind}\nolimits}
\def\add{\mathop{\rm add}\nolimits}

\newcommand{\DBH}{D^b(H)}

\newcommand{\Ga}{\Gamma}

\newcommand{\C}{\mathcal{C}}
\newcommand{\A}{\mathcal{A}}
\newcommand{\B}{\mathcal{B}}
\newcommand{\E}{\mathcal{E}}
\newcommand{\LL}{\mathcal{L}}
\newcommand{\RR}{\mathcal{R}}

\def\la{{\mathcal L}_A}
\def\ra{{\mathcal R}_A}

\newcommand{\mor}[3]{\xymatrix@1@C=15pt{#1: #2\ar[r]& #3}}
\newcommand{\path}[3]{\xymatrix@1@C=15pt{#1: #2\ar@{~>}[r]& #3}}

\begin{document}

\title{On tilting modules over cluster-tilted algebras}%

\author[D.~Smith]{David Smith}

\address{Department of Mathematical Sciences, Norwegian
University of Science and Technology, N-7491, Trondheim, Norway}
\email{unracinois@hotmail.com}

\thanks{Research supported by NSERC of Canada and by the Norwegian Research Council (Storforsk grant
no.~167130)}



\begin{abstract}
In this paper, we show that the tilting modules over a
cluster-tilted algebra $A$ lift to tilting objects in the associated
cluster category $\mathcal{C}_H$. As a first application, we
describe the induced exchange relation for tilting $A$-modules
arising from the exchange relation for tilting object in
$\mathcal{C}_H$. As a second application, we exhibit tilting
$A$-modules having cluster-tilted endomorphism algebras.
\end{abstract}

\maketitle

%
%

Cluster algebras were introduced by Fomin and Zelevinsky \cite{FZ02}
in the context of canonical basis of quantized enveloping algebras
and total positivity for algebraic groups, but quickly turned out to
be related to many other fields in mathematics. In the
representation theory of finite dimensional algebras, the so-called
cluster categories were introduced in \cite{BMRRT06} (and also in
\cite{CCS06} for the $\mathbb{A}_n$ case) as a natural categorical
model for the combinatorics of the corresponding cluster algebras of
Fomin and Zelevinsky. The construction is as follows. Let $Q$ be a
quiver without oriented cycles. There is then, for a field $k$, an
associated finite dimensional hereditary path algebra $H=kQ$. Since
$H$ has finite global dimension, its bounded derived category $\DBH$
of the finitely generated modules has almost split triangles
\cite{H88}. Let $\tau$ be the corresponding translation functor.
Denoting by $F$ the composition $\tau^{-1}[1]$, where $[1]$ is the
shift functor in $\DBH$, the cluster category $\C_H$ was defined as
the orbit category $\DBH/F$, and was shown to be canonically
triangulated \cite{K05} and have almost split triangles
\cite{BMRRT06}.

In this model, the exceptional objects are associated with the
cluster variables of \cite{FZ02} while the tilting objects
correspond to the clusters.
Remarkably, one also defines an exchange relation on the tilting
objects in $\C_H$, corresponding to the exchange relation on the
clusters of \cite{FZ02}.  More precisely, an almost complete tilting
object $\overline{T}$ in $\C_H$ has exactly two nonisomorphic
indecomposable complements $M$ and $M^*$, and these are related by
triangles
$$\xymatrix@1@C=15pt{M^* \ar[r]^g & B \ar[r]^f & M \ar[r] &
M^*[1]}
\mbox{\qquad and \qquad}%
\xymatrix@1@C=15pt{M \ar[r]^{f^*} & B^* \ar[r]^{g^*} & M^* \ar[r]
& M[1]}$$
where $f,g^*$ are minimal right $\add \overline{T}$-approximations
and $f^*,g$ are minimal left $\add \overline{T}$-approximations
(see \cite{BMRRT06}).

In view of the importance of tilting theory in the representation
theory of finite dimensional algebras, the (opposite) endomorphism
algebras of these tilting objects, called \textit{cluster-tilted},
were then introduced and studied in \cite{BMR07} (see also
\cite{CCS06}). Their module theory was shown to be to a large
extent determined by the cluster categories in which they arise.
Indeed, given a cluster category $\C_H$ and a tilting object $T$
in $\C_H$, it was shown by Buan, Marsh and Reiten \cite{BMR07}
that the functor $\Hom_{\C_H}(T,-)$ induces an equivalence
$\xymatrix@1@C=15pt{\C_H/{\add T[1]} \ar[r] &
\mod\End_{\C_H}(T)^{op}}$.

Since then, cluster-tilted algebras have been studied by several
authors, and revealed to have very nice properties, see for
instance \cite{ABS08, ABS07, BMR05, KR07}. In particular, they
were shown in \cite{KR07} to be Gorenstein of dimension at most
one and in \cite{ABS07} to be obtained from tilted algebras by
trivial extensions.

In this paper, we are interested in the problem of identifying
tilting modules over cluster-tilted algebras.  Our motivation
comes from two points of view.  On one side, the nice exchange
relation for tilting objects over cluster categories should carry
over Buan-Marsh-Reiten's equivalence and result in a similar
exchange relation for tilting modules over cluster-tilted
algebras, allowing to identify many tilting modules. Of course,
one then has to care about projective dimensions. On the other
hand, as stressed above, cluster-tilted algebras enjoy some very
nice properties.  Tilting theory being intimately related to
derived equivalences (under which many properties are known to be
preserved) by Happel's and Rickard's Theorems \cite{H88, R89}, the
study of tilting modules is then a natural question.

In what follows, we present two different methods to find tilting
modules over cluster-tilted algebras, dividing the paper in two
distinct parts.

The first approach follows the above discussion, in the sense that
we study the exchange relation of tilting modules over
cluster-tilted algebras coming from the exchange relation of
tilting objects for cluster categories. As pointed out above, one
then has to care about projective dimension in the following
sense: if $T$ and $T'$ are two tilting objects over a cluster
category $\C_H$ such that $\add T[1] \cap \add T'=\{0\}$, then it
follows from Buan-Marsh-Reiten's equivalence (see also \cite{KR07,
KZ07}) that the image of $T'$ under the equivalence
$\xymatrix@1@C=15pt{\C_H/{\add T[1]} \ar[r] &
\mod\End_{\C_H}(T)^{op}}$ is exceptional and has the right number
of indecomposable direct summands to be a tilting module, but a
priori no one knows about its projective dimension, which
generally turns out to be infinite. The situation is better in the
other direction. Indeed, while lifting tilting modules over
cluster-tilted algebras to objects in the cluster category
obviously does not bring any projective dimension problems, one
now has to care about the exceptionality of the resulting objects,
since the cluster category contains more maps, namely those
factoring through $\add T[1]$.  The following theorem says that
such problems do not occur.
We stress that by abuse of notation, we also denote, here and in
the sequel, by $M$ the preimage in $\C_H$ of an
$\End_{\C_H}(T)^{op}$-module $M$ under the composition
$\xymatrix@1@C=15pt{\C_H\ar[r]&\C_H/\add T[1] \ar[r] &
\mod\End_{\C_H}(T)^{op}}$.

\begin{my_thm}\label{Theorem 1}%
Let $\C_H$ be a cluster category, $T$ be a tilting object in
$\C_H$ and $A=\End_{\C_H}(T)^{op}$. Let $M, N$ be $A$-modules of
projective dimension at most one. If $\Ext^1_{A}(M,N)=0$ and
$\Ext^1_{A}(N,M)=0$, then $\Ext^1_{\C_H}(M, N)=0$ and
$\Ext^1_{\C_H}(N, M)=0$. In particular, the tilting $A$-modules
lift to tilting objects in $\C_H$.
\end{my_thm}

From this, we get that the endomorphism algebras of tilting
modules over cluster-tilted algebras are quotients of
cluster-tilted algebras (Corollary \ref{cor quotient}).

On the other hand, the study of the possible complements for an
almost complete tilting module has been the central point of many
investigations during the past years.  It is known that an almost
complete tilting module of projective dimension at most one admits
at most two nonisomorphic complements. Combining Theorem 1 with a
result from \cite{CHU94, H95} (see Theorem \ref{thm CHU}) then
allows to show that for a cluster-tilted algebra, these two
complements are related by the exchange relation in $\C_H$.

\begin{my_prop}\label{Proposition 2}Let $\C_H$ be a cluster category,
$T$ be a tilting object in $\C_H$ and $A=\End_{\C_H}(T)^{op}$. Let
$S=\overline{S}\oplus M$ be a (basic) tilting $A$-module, with $M$
indecomposable.
Also, let %
$$\xymatrix@1@C=15pt{M^* \ar[r]^g & B \ar[r]^f & M \ar[r] &
M^*[1]}
\mbox{\qquad and \qquad}%
\xymatrix@1@C=15pt{M \ar[r]^{f^*} & B^* \ar[r]^{g^*} & M^* \ar[r]
& M[1]}$$
be the corresponding exchange triangles in $\C_H$, where $f,g^*$
are minimal right $\add \overline{S}$-approximations in $\C_H$ and
$f^*,g$ are minimal left $\add \overline{S}$-approximations in
$\C_H$. The following are equivalent:
\begin{enumerate}
\item[(a)] There exists an indecomposable module $M'$, not isomorphic to $M$, such that $\overline{S}\oplus M'$ is a tilting
$A$-module;
\item[(b)] $\overline{S}\oplus M^*$ is a tilting $A$-module;
\item[(c)] As an $A$-module, $M^*\neq 0$ and $\pd_AM^*\leq 1$.
\item[(d)] Either $\Hom_{\C_H}(T,f)$ is an epimorphism in $\mod A$ or $\Hom_{\C_H}(T,f^*)$ is a monomorphism in
$\mod A$;
\item[(e)] $\overline{S}$ is a faithful $A$-module.
\end{enumerate}
\end{my_prop}

The second method deals with completely different tools. Given an
algebra $A$, we consider the left part $\la$ and the right part
$\ra$ of its module category $\mod A$ (see \cite{HRS96}). In
\cite{ACT04}, Assem, Coelho and Trepode studied the algebras $A$
for which the subcategory $\add\la$ is functorially finite in
$\mod A$ (in the sense of \cite{AS80}) and called them
\textit{left supported}. Dually, they defined the \textit{right
supported} algebras. They proved that $A$ is left supported if and
only if a specific $A$-module $L$ is a tilting module, and
similarly for the right supported algebras. As we shall see, the
left and the right parts of a cluster-tilted not hereditary
algebra are both finite, implying that any cluster-tilted algebra
is left and right supported. The module $L$ is the direct sum of
the indecomposable Ext-injective modules in $\add\la$ and the
indecomposable projective modules which are not in $\la$.  Hence
$L$ determines a "slice" in $\la$ given by the sum of the
indecomposable Ext-injective modules in $\add\la$. Our results
show that any basic object $S$ in $\add\la$, which is maximal for
the property that $\Ext^1_A(S,S)=0$, gives rise to a tilting
module. However, the ones given by slices in $\la$, called
\textit{$\la$-slices} (see Definition \ref{defn la-slice}) give
remarkable tilting modules, since their endomorphism algebra is
still cluster-tilted.

\begin{my_thm}\label{Theorem 3} %
Let $\C_H$ be a cluster category, $T$ be
a tilting object in $\C_H$ and $A$ be the cluster-tilted algebra
$\End_{\C_H}(T)^{op}$. Assume that $A$ is not hereditary and let
$\Sigma$ be an $\la$-slice. Also, let $F=\oplus_{i=1}^m P_i$
denote the direct sum of all indecomposable projective modules not
in $\la$. Then,
\begin{enumerate}
\item[(a)] $T_{\Sigma}=\Sigma\oplus F$ is a tilting
$A$-module;
\item[(b)] The algebra $A_\Sigma=\End_A(T_{\Sigma})^{op}$ is
isomorphic to $\End_{\C_H}(T_{\Sigma})^{op}$.  In particular,
$A_\Sigma$ is cluster-tilted;
\item[(c)] The quiver of $A_\Sigma$ is obtained from that of $A$ with a finite number of reflections at sinks.
\end{enumerate}
\end{my_thm}


This paper is organized as follows.  In Section 1, we collect the
necessary background concerning cluster categories and
cluster-tilted algebras. The Sections 2 and 3 are devoted to the
proofs of Theorem \ref{Theorem1} and Proposition \ref{Proposition
2} respectively. Finally, after some necessary preliminaries on
supported algebras in Section 4, we prove Theorem \ref{Theorem 3}
in Section 5.

This work was completed when the author was visiting NTNU in
Norway as a postdoc.  The author would like to thank I.~Reiten and
A.~Buan for some fruitful discussions.

%
%
\section{First preliminaries}
    \label{First preliminaries}

In this section, we review some useful notions and results that
will be used for the proofs of Theorem \ref{Theorem 1} and
Proposition \ref{Proposition 2}. More preliminaries concerning
Theorem \ref{Theorem 3} are postponed to Section \ref{More
preliminaries}.

%
%
\subsection{Basic notations}
    \label{Basic notations}

In this paper, all algebras are connected finite dimensional
algebras over a field $k$. For an algebra $A$, we denote by $\mod
A$ the category of finitely generated (right) $A$-modules. For an
$A$-module $M$, we respectively denote by $\pd_A M$ and $\id_AM$
the projective dimension and the injective dimension of $M$.

More generally, for an additive category $\A$ we let $\ind \A$ be
a full subcategory whose objects are representatives of the
isomorphism classes of indecomposable objects in $\A$. By an
indecomposable object in $\A$ we therefore mean an object in
$\ind\A$.  In case $\A=\mod A$, for some algebra $A$, we write
$\ind A$ instead of $\ind(\mod A)$. For an object $T$ in $\A$,
$\add T$ denotes the full subcategory of $\A$ with objects all
direct summands of direct sums of copies of $T$.

Also, even though the notions of tilting object slightly differ
according to the type of categories we consider (see Sections
\ref{Cluster categories and tilting objects} and
\ref{Cluster-tilted algebras and tilting modules} for details), we
will in any case say that an object $\overline{T}$ in $\A$ is an
\textit{almost complete tilting object} if it is not a tilting
object but there exists an indecomposable object $M$ in $\A$ such
that $\overline{T}\oplus M$ is a tilting object. In this case, $M$
is said to be a \textit{complement} for $\overline{T}$. Finally,
all (partial) tilting objects $T$ we consider are assumed to be
\textit{basic}, that is if $T=\oplus_{i=1}^nT_i$ is a
decomposition in indecomposable direct summands of $T$, then
$i\neq j$ implies $T_i\ncong T_j$.

%
%
\subsection{Approximations}
    \label{Approximations}

Let $\A$ be an additive category and $\B$ be an additive
subcategory of $\A$.  For an object $A$ in $\A$, a map
$\mor{f}{B}{A}$, with $B\in\B$ is called a \textit{right
$\B$-approximation} if any other map $\mor{f'}{B'}{A}$ with $B'\in
\B$ factors through $f$, that is there exists $\mor{g}{B'}{B}$
such that $f'=fg$. There is the dual notion of \textit{left
$\B$-approximation}. If any object in $\A$ admits a right (left)
$\B$-approximation, then $\B$ is said to be a
\textit{contravariantly (covariantly) finite} subcategory of $\A$.
It is called \textit{functorially finite} if it is both
contravariantly finite and covariantly finite. Finally a
\textit{minimal right $\B$-approximation} is a right
$\B$-approximation $\mor{f}{B}{A}$ such that for every
$\mor{g}{B}{B}$ such that $fg=f$, the map $g$ is an isomorphism.
The \textit{minimal left $\B$-approximation} are defined dually.
These notions were introduced in \cite{AS80}.

%
%
\subsection{Cluster categories and tilting objects}
    \label{Cluster categories and tilting objects}

Let $H$ be a hereditary algebra. As mentioned in the introduction,
the cluster category $\C_H$ is the orbit category $\DBH/F$, where
$F=\tau^{-1}[1]$. Thus, the objects in $\C_H$ are the $F$-orbits $X
= (F^i\widetilde{X})_{i\in\mathbb{Z}}$, where $\widetilde{X}$ is an
object in $\DBH$. The set of morphisms from $X =
(F^i\widetilde{X})_{i\in\mathbb{Z}}$ to $Y =
(F^i\widetilde{Y})_{i\in\mathbb{Z}}$ in $\C_H$ is given by
$$\Hom_{\C_H}(X,Y)=\prod_{i\in\mathbb{Z}}
\Hom_{\DBH}(\widetilde{X}, F^i\widetilde{Y}).$$
It is shown in \cite{K05} that $\C_H$ is a triangulated category
and that the canonical functor $\xymatrix@1@C=15pt{\DBH \ar[r]&
\C_H}$ is a triangle functor. Moreover, $\C_H$ has almost split
triangles and $\tau_{\C_H}=[1]$. Let $\mathcal{F}=\ind(\mod H \cup
H[1])$, that is the set consisting of the indecomposable
$H$-modules together with the objects $P[1]$ where $P$ is an
indecomposable projective $H$-module.  It is easily seen that
$\mathcal{F}$ contains exactly one representative from each
$F$-orbit of indecomposable objects in $\DBH$. Hence,
$\mathcal{F}=\ind\C_H$ and we can (and will) always assume that an
indecomposable object in $\C_H$ is a $H$-module or of the form
$P[1]$. Moreover, for two objects $M, N$ in $\mathcal{F}$, we have
$\Hom_{\DBH}(\widetilde{M},F^i\widetilde{N})=0$ for all $i\neq 0,
1$ (see \cite[(1.5)]{BMRRT06}). Also, by
\cite[(1.4)(1.7)]{BMRRT06},
$$D\Ext^1_{\C_H}(N,M)\cong \Ext^1_{\C_H}(M,N) \cong
D\Hom_{\C_H}(N, \tau M)$$

Let $T$ be a basic object in $\C_H$. Following \cite{BMRRT06}, $T$
is a \textit{cluster-tilting object}, or a \textit{tilting object}
for short, provided $\Ext^1_{\C_H}(T,T)=0$ and $T$ has a maximal
number of nonisomorphic direct summands (corresponding to the
number of nonisomorphic simple $H$-modules). Moreover, up to
derived equivalence, one can always assume that a given tilting
object $T$ is induced by a tilting module over $H$ (see
\cite[(3.3)]{BMRRT06}). Also,  an almost complete basic tilting
object $\overline{T}$ in $\C_H$ has exactly two nonisomorphic
complements $M$ and $M^*$, and these are related by some exchange
triangles
$$\xymatrix@1@C=15pt{M^* \ar[r]^g & B \ar[r]^f & M \ar[r] &
M^*[1]}
\mbox{\qquad and \qquad}%
\xymatrix@1@C=15pt{M \ar[r]^{f^*} & B^* \ar[r]^{g^*} & M^* \ar[r]
& M[1]}$$
where $f,g^*$ are minimal right $\add \overline{T}$-approximations
and $f^*,g$ are minimal left $\add \overline{T}$-approximations.
The following particular case will be heavily exploited in Section
5. For more details on cluster categories, we refer to
\cite{BMRRT06}.

\begin{my_rem}\label{rem ase}
Let $\overline{T}$, $M$ and $B$ be as above and let
$\xymatrix@1@C=15pt{M\ar[r]&Q\ar[r]&\tau^{-1}M\ar[r]&M[1]}$ be the
almost split triangle starting at $M$.  If $Q\in\add\overline{T}$,
then $Q=B$ and therefore $M^*=\tau^{-1}M$.  Hence the exchange of
$M$ by $M^*$ coincides with an almost split exchange in $\C_H$.
\end{my_rem}
%

%
%
\subsection{Cluster-tilted algebras and tilting modules}
    \label{Cluster-tilted algebras and tilting modules}

A \textit{cluster-tilted} algebra is an algebra of the form
$A=\End_{\C_H}(T)^{op}$, for some tilting object $T$ in a cluster
category $\C_H$.  Moreover, by \cite{BMR07}, the functor
$\Hom_{\C_H}(T,-)$ induces an equivalence
$\xymatrix@1@C=15pt{\C_H/\add T[1] \ar[r] & \mod A}$ under which
the almost split sequences in $\mod A$ are induced by almost split
triangles in $\C_H$. Moreover, it was shown in \cite{KR07} that
any cluster-tilted algebra $A$ is Gorenstein of dimension at most
one, meaning that every projective module is of injective
dimension at most one, and dually every injective module is of
projective dimension at most one. As an important consequence, the
projective dimension and the injective dimension of any $A$-module
are simultaneously either infinite, or less or equal than one (see
\cite[(Section 2.1)]{KR07}). In particular, the tilting modules
are of projective dimension at most one. Hence a (basic)
$A$-module $S$ is a \textit{tilting $A$-module} if :
\begin{itemize}
\item $\pd_{A}S\leq 1$ (equivalently $\id_{A}S\leq 1$);
\item $\Ext_{A}^1(S,S)=0$;
\item The number of indecomposable direct summands of $S$ equals the number of simple
$A$-modules (equivalently simple $H$-modules).
\end{itemize}
Also, we recall that in this paper, we keep the same notation for
an $A$-module and its preimage in $\C_H$ under the projection
$\xymatrix@1@C=15pt{\C_H\ar[r]&\C_H/\add T[1] \ar[r] & \mod A}$.

%
%
\section{Proof of Theorem \ref{Theorem 1}}
    \label{Proof of Theorem 1}

In this section, we recall some useful features of modules of
projective or injective dimension at most one and prove Theorem
\ref{Theorem 1}. We start with the following well-known lemma (see
\cite[(IV.2.13)(IV.2.14)]{ASS06} for instance).

\begin{my_lem}\label{lem pd1}%
Let $A$ be an algebra and $M$ be an
$A$-module. %
\begin{enumerate}
\item[(a)] $\pd_AM \leq 1$ if and only if $\Hom_A(DA, \tau M)=0$.
Moreover, if $\pd_AM\leq 1$, then $\Ext^1_A(M,N)\cong D\Hom_A(N,\tau
M)$ for each $A$-module $N$;
\item[(b)] $\id_AM \leq 1$ if and only if $\Hom_A(\tau^{-1}M, A)=0$.
Moreover, if $\id_AM\leq 1$, then $\Ext^1_A(N,M)\cong
D\Hom_A(\tau^{-1}M,N)$ for each $A$-module $N$;
\end{enumerate}
where $\mor{D=\Hom_k(-,k)}{\mod A^{op}}{\mod A}$ denotes the usual
duality.
\end{my_lem}

We note that if $\C_H$ is a cluster category and $T$ is a tilting
object in $\C_H$, with $A=\End_{\C_H}(T)^{op}$, then the
equivalence $\xymatrix@1@C=15pt{\C_H/\add T[1] \ar[r] & \mod A}$
takes the objects in $\add T$ to projective $A$-modules and the
objects in $\add T[2]$ to injective $A$-modules. In view of this
and the Gorenstein property of cluster-tilted algebras, the above
lemma immediately implies the following result.

\begin{my_lem}\label{lem pd2}%
Let $\C_H$ be a cluster category and $T$ be
a tilting object in $\C_H$.  Let $A=\End_{\C_H}(T)^{op}$ and $M$
be an $A$-module.  The following conditions are equivalent:%
\begin{enumerate}
\item[(a)] $\pd_{A}M \leq 1$;
\item[(b)] In $\C_H$, any map from an object in $\add T[2]$ to $M[1]$ factors through $\add T[1]$;
\item[(c)] $\id_{A}M \leq 1$;
\item[(d)] In $\C_H$, any map from $M[-1]$ to an object in $\add T$ factors through $\add
T[1]$.
\end{enumerate}
\end{my_lem}

We are now in position to prove Theorem \ref{Theorem 1}.

\begin{my_thm2}\label{Theorem1}%
Let $\C_H$ be a cluster category, $T$ be a tilting object in
$\C_H$ and $A=\End_{\C_H}(T)^{op}$. Let $M, N$ be $A$-modules of
projective dimension at most one. If $\Ext^1_{A}(M,N)=0$ and
$\Ext^1_{A}(N,M)=0$, then $\Ext^1_{\C_H}(M, N)=0$ and
$\Ext^1_{\C_H}(N, M)=0$. In particular, the tilting $A$-modules
lift to tilting objects in $\C_H$.
\end{my_thm2}
\begin{proof}
Clearly, it suffices to prove the Theorem for $M,N$
indecomposable. Moreover, we assume that $T$ is induced by a
tilting $H$-module.

We first assume that, in $\C_H$, $M$ and $N$ are two $H$-modules.
Also, assume to the contrary that $\Ext_{\C_H}^1(M,N)\neq 0$.
Then,
\begin{eqnarray*}
0\neq \Ext_{\C_H}^1(M,N) & \cong & D\Hom_{\C_H}(N,\tau M) \\
&=& D\Hom_{\DBH}(N,\tau M) \oplus D\Hom_{\DBH}(N,M[1])\\
&\cong& \Hom_{\DBH}(M,N[1])\oplus D\Hom_{\DBH}(N,M[1])
\end{eqnarray*}
and thus $\Hom_{\DBH}(M,N[1])\neq 0$ or $\Hom_{\DBH}(N,M[1])\neq 0$.
Assume, without loss of generality, that $\Hom_{\DBH}(M,N[1])\neq
0$.  Also, we have
\begin{eqnarray*}
0= D\Ext_{A}^1(N,M) & \cong & \Hom_{A}(M,\tau N) \\
&\cong & \frac{\Hom_{\C_H}(M,N[1])}{\{\mor{f}{M}{N[1]} \mbox{
factoring through }\add T[1]\}}
\end{eqnarray*}
where the first isomorphism follows from Lemma \ref{lem pd1}.
Therefore, any map in $\Hom_{\C_H}(M,N[1])$ factors through $\add
T[1]$. Similarly, any map in $\Hom_{\C_H}(N,M[1])$ factors through
$\add T[1]$. Lifting this property to $\DBH$ means, in particular,
that any map in $\Hom_{\DBH}(M,N[1])$ factors through $\add(\tau
T\oplus T[1])$. Now, let $\{f_1,\dots, f_n\}$ be a basis for
$\Hom_{\DBH}(M,N[1])$. For each $i$, there exist $T'_i, T''_i$ in
$\add T$ and maps $\mor{(\substack{\alpha_i\\ \beta_i})}{M}{\tau
T'_i \oplus T''_i[1]}$ and $\mor{(\gamma_i, \delta_i)}{\tau T'_i
\oplus T''_i[1]}{N[1]}$ such that $f_i=(\gamma_i, \delta_i)\circ
(\substack{\alpha_i\\ \beta_i})$. Taking $T'=\oplus_{i=1}^{n}T'_i$,
$T''=\oplus_{i=1}^{n}T''_i$, $\alpha=\mbox{diag}(\alpha_1, \dots,
\alpha_n)$ and $\beta=\mbox{diag}(\beta_1, \dots, \beta_n)$, we see
that any map
in $\Hom_{\DBH}(M,N[1])$ factors through $\mor{(\substack{\alpha\\
\beta})}{M}{\tau T'\oplus T''[1]}$. In other words, we have a
surjective map
$$\xymatrix@1@C=30pt{\Hom_{\DBH}(\tau T', N[1])\oplus
\Hom_{\DBH}(T''[1], N[1]) \ar[r]^-{\circ (\substack{\alpha\\
\beta})} & \Hom_{\DBH}(M, N[1])}$$
Under the natural isomorphism
$$\Hom_{\DBH}(X,Y[1])\cong \Ext_H^{1}(X,Y)\cong D\Hom_H(Y, \tau
X),\mbox{ for }X,Y\in\mod H, $$
the map $\mor{\beta}{M}{T''[1]}$ becomes an element of
$D\Hom_H(T'', \tau M)$. More generally the above surjective map
becomes the surjective map
$$\xymatrix@1@C=30pt{D\Hom_{H}(N, \tau^2 T')\oplus \Hom_{H}(T'',
N) \ar[r] & D\Hom_{H}(N, \tau M)}$$
which takes a pair $(g,h)$ in $D\Hom_{H}(N, \tau^2T')\oplus
\Hom_{H}(T'', N)$ to the morphism $\xymatrix@1@C=15pt{\Hom_H(N,\tau
M)\ar[r] & k}$ sending an element $f\in \Hom_A(N,\tau M)$ onto the
element $g(\tau(\alpha)f)+\beta(fh)$.
Applying the duality $D$ yields an injective map
$$\xymatrix@1@C=30pt{\Hom_{H}(N, \tau M)\ar[r] & \Hom_{H}(N,
\tau^2 T')\oplus D\Hom_{H}(T'', N)}$$
taking an element $g\in \Hom_H(N,\tau M)$ to the pair
$(\tau(\alpha)g, \overline{g})$, where $\overline{g}(h)=\beta(gh)$
for $h\in \Hom_H(T'',N)$.
Now, recall that by assumption $0\neq \Hom_{\DBH}(M, N[1])\cong
\Hom_{H}(N, \tau M)$. Hence, let $g$ be a nonzero morphism in
$\Hom_{H}(N, \tau M)$.  The injectivity of the above map gives
$\tau(\alpha)g\neq 0$ or $gh\neq 0$ for some $h\in \Hom_H(T'',N)$.
In other words, one of the two compositions
$$\xymatrix@1@C=25pt{N\ar[r]^-{g}&\tau M \ar[r]^-{\tau(\alpha)} &
\tau^{2}T'}
\mbox{\qquad and \qquad}%
\xymatrix@1@C=25pt{T''\ar[r]^-{h} & N \ar[r]^-{g} & \tau M}$$
is not zero. However, since any map in
$\Hom_{\C_H}(N,M[1])=\Hom_{\DBH}(N, \tau M)\oplus
\Hom_{\DBH}(N,M[1])$ factors through $\add T[1]$, $g$ factors
through $\add \tau T$ in $\mod H$, say through $\tau T'''$, with
$T'''\in\add T$. The above compositions then yield a nonzero map
of the form $\xymatrix@1@C=15pt{\tau T'''\ar[r] & \tau^{2}T'}$ or
$\xymatrix@1@C=15pt{T''\ar[r] & \tau T'''}$, a contradiction to
$\Ext^1_{\C_H}(T,T)=0$. Hence $\Ext^1_{\C_H}(M,N)=0$, and dually
$\Ext^1_{\C_H}(N,M)=0$.

We now assume that $M\in\mod H$ and $N\in\add H[1]$. Let $P$ be an
indecomposable projective $H$-module such that $N=P[1]$. Then
$\tau N =I$, where $I$ is the indecomposable injective $H$-module
satisfying $\soc I = \top P$. Now, assume that
$\Ext^1_{\C_H}(M,N)\neq 0 \neq \Ext^1_{\C_H}(N,M)$. We have
$$0\neq
\Ext^1_{\C_H}(M,N)=\Ext^1_{\C_H}(M,P[1])=D\Hom_{\DBH}(P[1],M[1])$$
and so $\Hom_H(P,M)\neq 0$.  Similarly, $\Ext^1_{\C_H}(N,M)$
yields $\Hom_H(M,I)\neq 0$. Let $f \in \Hom_H(P,M)$ and $g\in
\Hom_H(M,I)$ be nonzero morphisms.  Since $\soc I = \top P$, we
get a nonzero composition
$\xymatrix@1@C=20pt{P\ar[r]^f&M\ar[r]^g&I}.$
Since $\pd_AM\leq 1$ and $\pd_AN\leq 1$, it follows from the first
part of the proof that $f$ factors through $\add T$ while $g$
factors through $\add\tau T$, contradicting $\Ext^1_H(T,T)=0$.
Hence $\Ext^1_{\C_H}(M,N)=0= \Ext^1_{\C_H}(N,M)$.
Finally, if $M$ and $N$ are both in $\add H[1]$, then
$\Ext^1_{\C_H}(M,N)=0= \Ext^1_{\C_H}(N,M)$ and we are done.
\end{proof}

The following easy example shows that Theorem \ref{Theorem 1} is
no longer true if we drop the assumption that $\pd_{A}M\leq 1$ and
$\pd_{A}N\leq 1$.

\begin{my_ex}\label{ex nolonger}%
Let $Q$ be the quiver $\xymatrix@C=10pt{1\ar[r]&2\ar[r]&3}$ and
$H=kQ$ be the path algebra.  The AR-quiver of the corresponding
cluster category $\C_H$ is given by
$$\xymatrix@R=10pt@!C=5pt{ & & {\substack{1 \\2
\\ 3}} \ar[dr] & & 3[1] \ar[dr] & & 2[1] \cong 3
\ar[dr] & &  \\ & {\substack{2 \\ 3}} \ar[dr] \ar[ur] & &
{\substack{1 \\2}} \ar[dr] \ar[ur] & & {\substack{2
\\ 3}}[1] \ar[dr] \ar[ur] & & {\substack{1
\\2}}[1] \cong  {\substack{2 \\ 3}} \ar[dr] &  \\ 3 \ar[ur] & & 2 \ar[ur] & & 1
\ar[ur] & & {\substack{1 \\2 \\ 3}}[1] \ar[ur] & & 3[2] \cong
{\substack{1 \\2 \\ 3}}}$$%
Let $T$ be the tilting object $T=3\oplus \substack{1\\ 2\\
3}\oplus 1$ and $A=\End_{\C_H}(T)^{op}$ be the corresponding
(self-injective) cluster-tilted algebra. Its AR-quiver is given
by%
$$ \xymatrix@R=10pt@!C=5pt{ & & \fbox{$\substack{2 \\1}$} \ar[dr]
&&&& \fbox{$\substack{1\\ 3}$} \ar[dr]\\ \ar@{--}[r] & 1 \ar[ur]
\ar@{--}[rr] & & 2 \ar[dr] \ar@{--}[rr] & & 3 \ar[ur] \ar@{--}[rr]
&& 1 \ar[dr] \ar@{--}[r] &
\\ \fbox{$\substack{1\\ 3}$} \ar[ur] &&&& \fbox{$\substack{3\\ 2}$} \ar[ur]&&&& \fbox{$\substack{2\\
1}$}}$$
where the dashed lines represent the Auslander-Reiten translates
and the identified modules are the projective-injective modules.
In $\mod A$, let $M=\substack{1\\ 3}$ and $N= 2$. Since $M$ is
projective-injective, we have $\Ext_A^1(M,N)=\Ext_A^1(N,M)=0$. In
$\C_H$, $M$ and $N$ correspond to the objects $3$ and
$\substack{1\\ 2}$ respectively, and
$\Ext^1_{\C_H}(M,N)=\Hom_{\C_H}(M, N[1])\neq 0$. Since $\pd_AM\leq
1$ but $\pd_AN=\infty$, this shows that Theorem \ref{Theorem1} is
no longer true when we drop the assumption that $\pd_{A}M\leq 1$
and $\pd_{A}N\leq 1$. Also, consider $M=\substack{1\\3}$ and
$N'=1$ in $\mod A$. Then $\pd_AN'=\infty$ but
$\Ext^1_{\C_H}(M,N')=\Ext^1_{\C_H}(N',M)=0$ in $\C_H$, showing
also that the converse of Theorem \ref{Theorem1} generally fails.
\end{my_ex}

As a direct consequence of Theorem \ref{Theorem1}, we obtain the
following nice result.

\begin{my_cor}\label{cor quotient}%
Let $\C_H$ be a cluster category, $T$ be a tilting object in
$\C_H$ and $A=\End_{\C_H}(T)^{op}$. Let $S$ be a tilting
$A$-module.  Then $\End_{A}(S)^{op}$ is a quotient of the
cluster-tilted algebra $\End_{\C_H}(S)^{op}$.
\end{my_cor}
\begin{proof}
By Theorem \ref{Theorem 1}, $S$ is a tilting object in $\C_H$.
Hence $\End_{\C_H}(S)^{op}$ is cluster-tilted. The result then
from the equivalence $\C_H/\add \tau T \cong \mod A$.
\end{proof}

In Section \ref{Special tilting modules}, we discuss examples
where $\End_{A}(S)^{op}\cong \End_{\C_H}(S)^{op}$.

%
%
\section{Exchange relation for cluster-tilted algebras}
    \label{Exchange relation for cluster-tilted algebras}

Here we discuss the induced exchange relation of tilting modules
over cluster-tilted algebras in view of Theorem \ref{Theorem 1}
and the exchange relation for tilting objects in the cluster
categories. For clear reasons (for instance when a cluster-tilted
algebra has projective-injective modules), it is not always
possible to exchange an indecomposable direct summand $M$ of a
tilting module $\overline{S}\oplus M$ by another indecomposable
$M^*$ such that $\overline{S}\oplus M^*$ is a tilting module.  In
this section, we give sufficient and necessary conditions for the
existence of such a complement $M^*$ for cluster-tilted algebras.
Basically, we show that if such a $M^*$ exists, then it is given
by the exchange triangles in $\C_H$.

More generally, complements of almost complete tilting modules (of
arbitrary finite projective dimension) over artin algebras have
been studied by several authors, in particular by Coelho, Happel
and Unger (see \cite{CHU94, H95} for instance). A very weak
version of one of their main results, but sufficient for our
purpose, goes as follows:

\begin{theorem}\label{thm CHU}%
\cite{CHU94, H95} Let $A$ be an artin algebra with finite
finitistic dimension. Let $\overline{S}$ be an almost complete
tilting module with $\pd_AM\leq 1$.
\begin{enumerate}
\item[(a)] If $\overline{S}$ is not faithful, then $\overline{S}$
admits a unique complement.
\item[(b)] If $\overline{S}$ is faithful, then $\overline{S}$
admits exactly two complements $M$ and $M'$ and there exists a
short exact sequence
$$\xymatrix@1@C=15pt{0\ar[r]&M\ar[r]^-{f} & C \ar[r]^-{g} & M'
\ar[r] & 0}$$
where $f$ is a minimal left $\add\overline{S}$-approximation and
$g$ is a minimal right $\add\overline{S}$-approximation.
\end{enumerate}
\end{theorem}

Below, we show that Proposition 2 is obtained by combining the
above Theorem with Theorem \ref{Theorem 1}. We need to recall one
further result, borrowed from \cite[(2.3)]{KZ07}.

\begin{my_lem}\label{KZ(2.3)}%
Let $\C_H$ be a cluster category, $T$ be a
tilting object in $\C_H$ and $A=\End_{\C_H}(T)^{op}$. Let
$\xymatrix@1@C=15pt{L\ar[r]^-{h} & M \ar[r]^-{f} & N \ar[r]^-{g} &
L[1]}$ be a triangle in $\C_H$. Then, in $\mod A$,
\begin{enumerate}
\item[(a)] $\Hom_{\C_H}(T,f)$ is a monomorphism if
and only if $\Hom_{\C_H}(T,h)=0$.
\item[(b)] $\Hom_{\C_H}(T,f)$ is an epimorphism if
and only if $\Hom_{\C_H}(T,g)=0$.
\end{enumerate}
\end{my_lem}

We are now able to prove Proposition 2. We mention that the
existence of the exchange triangles in the statement follows from
Theorem \ref{Theorem 1}.

\begin{my_prop2}\label{Proposition2}%
Let $\C_H$ be a cluster category,
$T$ be a tilting object in $\C_H$ and $A=\End_{\C_H}(T)^{op}$. Let
$S=\overline{S}\oplus M$ be a (basic) tilting $A$-module, with $M$
indecomposable.
Also, let %
$$\xymatrix@1@C=15pt{M^* \ar[r]^-{g} & B \ar[r]^-{f} & M
\ar[r]^-{h} & M^*[1]}
\mbox{\qquad and \qquad}%
\xymatrix@1@C=15pt{M \ar[r]^-{f^*} & B^* \ar[r]^-{g^*} & M^*
\ar[r]^-{h^*} & M[1]}$$
be the corresponding exchange triangles in $\C_H$, where $f,g^*$
are minimal right $\add \overline{S}$-approximations in $\C_H$ and
$f^*,g$ are minimal left $\add \overline{S}$-approximations in
$\C_H$. The following are equivalent:
\begin{enumerate}
\item[(a)] There exists an indecomposable module $M'$, not isomorphic to $M$, such that $\overline{S}\oplus M'$ is a tilting
$A$-module;
\item[(b)] $\overline{S}\oplus M^*$ is a tilting $A$-module;
\item[(c)] As an $A$-module, $M^*\neq 0$ and $\pd_AM^*\leq 1$.
\item[(d)] Either $\Hom_{\C_H}(T,f)$ is an epimorphism in $\mod A$ or $\Hom_{\C_H}(T,f^*)$ is a monomorphism in
$\mod A$;
\item[(e)] $\overline{S}$ is a faithful $A$-module.
\end{enumerate}
\end{my_prop2}
\begin{proof}
Clearly, the equivalence of (a) and (e) follows from Theorem
\ref{thm CHU}. The same theorem also shows that (b) implies (e),
while trivially (b) implies (c).

We now show that (c) implies (b) and (d). By the exchange relation
in $\C_H$, we know that $\overline{S}\oplus M^*$ is a tilting
object in $\C_H$.  In particular,
$\Ext^1_{\C_H}(\overline{S}\oplus M^*, \overline{S}\oplus M^*)=0$,
and so $\Ext^1_{A}(\overline{S}\oplus M^*, \overline{S}\oplus
M^*)=0$ (see \cite[(4.9)]{KZ07}).  Since $\pd_{A}M^*\leq 1$ by
assumption, $\overline{S}\oplus M^*$ is a tilting $A$-module. This
shows (b).  Now, by Theorem \ref{thm CHU}, there exists a short
exact sequence of the form
$$\xymatrix@1@C=15pt{0\ar[r] & M^* \ar[r] & C
\ar[r]^{\underline{j}} & M \ar[r] & 0}
\mbox{\qquad or \qquad}%
\xymatrix@1@C=15pt{0\ar[r] & M \ar[r]^{\underline{j}^*} & C^*
\ar[r] & M^* \ar[r] & 0}$$
where $C, C^*\in\add\overline{S}$. Assume that the first exact
sequence exists, and let $\mor{j}{C}{M}$ be a morphism in $\C_H$
such that $\Hom_{\C_H}(T,j)=\underline{j}$. Now, since
$\mor{f}{B}{M}$ is a right $\add\overline{S}$-approximation, there
exists $\mor{f'}{C}{B}$ such that $j=ff'$. Then,
$\underline{j}=\Hom_{\C_H}(T,j)=\Hom_{\C_H}(T,f)\circ
\Hom_{\C_H}(T,f')$, showing that $\Hom_{\C_H}(T,f)$ is an
epimorphism. Similarly, if the second short exact sequence exists,
then $\Hom_{\C_H}(T,f^*)$ is a monomorphism. This shows (d).

Conversely, (d) implies (c).  Indeed, assume for instance that
$\Hom_{\C_H}(T,f)$ is an epimorphism.  By Lemma \ref{KZ(2.3)}, we
have $\Hom_{\C_H}(T,h)=0$.  Hence $h[-1]$ factors through $\add
T$. Since $\pd_{A}M\leq 1$, Lemma \ref{lem pd2} (d) implies that
$h[-1]$ factors through $\add T[1]$. Thus, by \cite[(3.4)]{KZ07},
we get a short exact sequence
$$\xymatrix@1@C=15pt{0 \ar[r] & M^* \ar[r]^-{\underline{g}} & B
\ar[r]^-{\underline{f}} & M \ar[r] & 0}$$
in $\mod A$, where $\underline{f}=\Hom_{\C_H}(T,f)$ and
$\underline{g}=\Hom_{\C_H}(T,g)$. Since $\pd_{A}M\leq 1$ and
$\pd_{A}B\leq 1$, we get $\pd_{A}M^*<\infty$, and so
$\pd_{A}M^*\leq 1$. Moreover, $M^*\neq 0$ since
$M\notin\add\overline{S}$. Hence (d) implies (c).

Finally, we show that (e) implies (b).  Assume that $\overline{S}$
is faithful.  By Theorem \ref{thm CHU}, there exists an
indecomposable module $M'$, not isomorphic to $M$, such that
$\overline{S}\oplus M'$ is a tilting module.  By Theorem
\ref{Theorem 1}, $\overline{S}\oplus M'$ is a tilting object in
$\C_H$, and since $M'\neq M$, we infer that $M'=M^*$.
\end{proof}

%
%
\section{More preliminaries}
    \label{More preliminaries}

Here starts the second part of the paper, whose objective is to
exhibit some tilting modules over cluster-tilted algebras whose
endomorphism algebras are again cluster-tilted. This is achieved
with the help of Theorem \ref{Theorem 1}, but also with the
property of cluster-tilted algebras of being left and right
supported.  Here, we gather the necessary terminology for the rest
of the paper.

%
%
\subsection{Paths and cycles}
    \label{Paths and cycles}

Let $A$ be an algebra. A \textit{path in $\ind
A$}, or simply a \textit{path}, is a sequence %
$ \delta : \xymatrix@C=15pt{M=M_0 \ar[r]^-{f_1} & M_1
\ar[r]^-{f_2} & \cdots \ar[r]^-{f_t} & M_t = N}$
$(t\geq 0)$ where $M_i \in \ind A$ and $f_i$ is a nonzero morphism
for each $i$. In this case, we write
$\xymatrix@1@C=15pt{M \ar@{~>}[r] & N}$
and we say that $M$ is a \textit{predecessor} of $N$ and $N$ is a
\textit{successor} of $M$. If each $f_i$ is irreducible, then
$\delta$ is \textit{sectional} if it contains no triple $(M_{i-1},
M_i, M_{i+1})$ such that $\tau_A M_{i+1} = M_{i-1}$. A
\textit{refinement} of $\delta$ is a path
$\xymatrix@C=15pt{M=M'_0 \ar[r] & M'_1 \ar[r] & \cdots \ar[r] &
M'_s = N,}$
with $s\geq t$, with an injective order-preserving function
$\mor{\sigma}{\{1,\dots,t-1\}}{\{1,\dots, s-1\}}$ such that
$M_i=M'_{\sigma(i)}$ when $1\leq i\leq t-1$. Finally, a path
$\delta$ is a \textit{cycle} if $M=N$ and at least one $f_i$ is
not an isomorphism. A subquiver $\Sigma$ of a connected component
$\Ga$ of the AR-quiver of $A$ is \textit{acyclic} if it contains
no cycle and \textit{convex} if any path in $\Ga$ starting and
ending in modules in $\Sigma$ consists only of modules in
$\Sigma$.

%
%
\subsection{The left and right parts of a module category}
    \label{The left and right parts of a module category}

For an algebra $A$, we define the left part $\la$ and the right
part $\ra$ of $\mod A$ as follows (see \cite{HRS96}):
    $$\begin{array}{c}
      \la=\{\ M\in \mbox{ind}A:\mbox{pd}_A N\leq 1 \mbox{ for each
predecessor } N \mbox{ of } M\ \}, \\
      \ra=\{\ M\in \mbox{ind}A:\mbox{id}_A N\leq 1 \mbox{ for each
successor } N \mbox{ of } M\ \}. \
    \end{array}$$
Clearly, $\la$ is closed under predecessors and $\ra$ is closed
under successors. The left and the right parts have been used in
recent years to describe many classes of algebras, amongst them
the quasitilted and the laura algebras (see \cite{AC03}). The next
result is helpful to detect the modules lying in these parts.

\begin{my_lem}\label{AC03 (1.6)}%
\cite[(1.6)]{AC03} Let $A$ be an algebra.%
\begin{enumerate}
\item[(a)] $\la$
consists of the modules $M\in\ind A$ such that, if there exists a
path from an indecomposable injective module to $M$, then this
path can be refined to a path of irreducible morphisms, and any
such refinement is sectional.
\item[(b)]  $\ra$ consists of the modules $N\in\ind A$ such
that, if there exists a path from $N$ to an indecomposable
projective module, then this path can be refined to a path of
irreducible morphisms, and any such refinement is sectional.
\end{enumerate}
\end{my_lem}

%
%
\subsection{Left and right supported algebras}
    \label{Left and right supported algebras}

In \cite{ACT04}, Assem, Coelho and Trepode defined the
\textit{left (right) supported} algebras as the algebras $A$ for
which the additive subcategory $\add \la$ ($\add\ra$) is
functorially finite in $\mod A$ (see Section
\ref{Approximations}).
Trivially, any hereditary algebra is left and right
supported.
In what follows, we mainly focus on left supported algebras,
instead of right supported algebras, and leave the primal-dual
translation to the reader.

When dealing with left supported algebras, the Ext-injective
modules in $\add\la$ play a prominent role since they determine if
the algebra is left supported or not. Recall that a module
$M\in\la$ is \textit{$\Ext$-injective in $\add\la$} if
$\Ext^1_A(N,M)=0$ for each $N\in \la$, or equivalently if
$\tau^{-1}M\notin\la$. Then, by \cite[(3.1)]{ACT04}, the class
$\mathcal{E}$ of indecomposable Ext-injective modules in $\la$ is
the union of two disjoint subclasses: %
$$\begin{array}{rl} \E_1 = & \{M \in \la : \mbox{ there exists an
injective $I$ in $\ind A$ and a path } \xymatrix@1@C=15pt{I
\ar@{~>}[r] & M}\}\\
\E_2 = &\{M\in\la\backslash\E_1 :\mbox{ there exists a projective
$P\in\ind A\backslash\la$ and a }\\ & \mbox{sectional path }
\xymatrix@1@C=15pt{P \ar@{~>}[r] & \tau^{-1}M}\}
\end{array}$$
Hence $\E=\E_1\cup\E_2$ and we denote by $E$ (or $E_1$, or $E_2$)
the direct sum of all indecomposable $A$-modules lying in $\E$ (or
$\E_1$, or $\E_2$ respectively). We also denote by $F$ the direct
sum of a full set of representatives of the isomorphism classes of
indecomposable projective $A$-modules not lying in $\la$. We set
$L=E\oplus F$ and $U=E_1\oplus \tau^{-1}E_2\oplus F$. With these
notations, we have the following reformulation of
\cite[(3.3)(4.2)]{ACT04} and \cite[(5.4)]{ACPT07}.

\begin{theorem}\label{thm ACT}%
An algebra $A$ is left supported if and only if $L$ is a tilting
$A$-module, and this occurs if and only if $U$ is a tilting
$A$-module.
\end{theorem}

As we will see, any cluster-tilted algebra is left supported, and
so the above provides canonical tilting modules, whose
endomorphism algebras will turn out to be again cluster-tilted.
For instance, in the easiest (but unfortunately degenerate and not
interesting) case where $\la=\emptyset$, we get the trivial
tilting module $L=U=A$, whose endomorphism algebra is obviously
cluster-tilted. Hopefully, we often get $\la\neq \emptyset$. In
fact, it is easily verified that for an algebra $A$, we have
$\la\neq\emptyset$ if and only if the ordinary quiver of $A$ has a
sink.

%
%
\section{Special tilting modules}
    \label{Special tilting modules}

In this section, we prove Theorem \ref{Theorem 3}.  This is made
in several steps.  We start by proving that any cluster-tilted is
left (and right) supported.

%
%
\subsection{Cluster-tilted algebras are left supported}
    \label{Cluster-tilted algebras are left supported}

Let $A$ be a cluster-tilted algebra.  If $A$ is hereditary, then
$\add\la=\mod A$, and so $A$ is trivially left (and dually right)
supported. Our first aim is to show that this property still holds
for cluster-tilted not hereditary algebras.  We need the following
lemma.

\begin{my_lem}\label{lem nsr}%
Let $\C_H$ be a cluster category, $T$ be a tilting object in
$\C_H$ and $A=\End_{\C_H}(T)^{op}$. Assume that $A$ is not
hereditary. Then any connected component of the AR-quiver of $A$
either contains no projective modules and no injective modules, or
contains both projective modules and injective modules.
\end{my_lem}
\begin{proof}
Let $P$ be an indecomposable projective $A$-module. Let $\Ga_A$
denote the AR-quiver of $A$ and $\Ga$ be the connected component
of $\Ga_{A}$ containing $P$. Also, let $\Sigma$ be the maximal
full, connected and convex subquiver of $\Ga$ containing only
indecomposable projective modules, including $P$. Since $A$ is not
hereditary, then $\Sigma$ has less vertices then the number of
$\tau$-orbits in $\Ga$. Hence, there exists $P'$ in $\Sigma$
together with an irreducible morphism $\xymatrix@1@C=15pt{M\ar[r]
& P'}$ in $\Ga$, where $M$ is indecomposable not projective. By
construction, $M$ belongs to $\Ga$. Moreover, let $T'$ be the
indecomposable direct summand of $T$ corresponding to $P'$. Since
$M$ is not projective, there is, in $\C_H$, an irreducible
morphism from $\tau^{2}T'$ to (the preimage of) $\tau M$. But
then, in $\Ga$, this corresponds to an irreducible morphism from
an indecomposable injective $A$-module $I$ to $\tau M$.  Hence
$\Ga$ contains at least one injective module. Dually, any
connected component containing an injective module also contains a
projective module.
\end{proof}

\begin{proposition}\label{prop ls}%
Let $A$ be a cluster-tilted not hereditary algebra.  Then $\la$
and $\ra$ are finite sets. In particular, $A$ is left and right
supported.
\end{proposition}
\begin{proof}
Assume that $\la\neq \emptyset$.  Since $\la$ is closed under
predecessors, $\la$ contains projective modules.  Let $P$ be such
a module. By \cite[(1.1)]{CL02} and Lemma \ref{lem nsr}, there
exists an integer $m\geq 0$ such that $\tau^{-m}P$ is a successor
of an injective module. Let $m$ be minimal for this property.
Then, by Lemma \ref{AC03 (1.6)}, we have $\tau^{-m-1}P\notin \la$
and so $\tau^{-m}P\in\E$. Since this holds for any projective in
$\la$, this shows that $A$ is left supported by
\cite[(3.3)]{ACT04}, and that $\la$ is finite by
\cite[(5.4)]{ACT04}. Dually, $\ra$ is finite.
\end{proof}

As a consequence, we get a straightforward characterization of the
cluster-tilted algebras which are laura. Recall from \cite{AC03}
that an algebra $A$ is \textit{laura} provided the set $\ind
A\setminus(\la\cup\ra)$ is finite. Therefore, a cluster-tilted
algebra is laura if and only if it is hereditary or representation
finite.

\begin{my_ex}\label{ex1}%
Let $A$ be the cluster-tilted algebra (of type $\mathbb{A}_8$)
given by the quiver
$$\xymatrix@1@C=8pt@R=8pt{&&&&&&&&& \bullet \ar[dr]^{\beta}\\ %
\bullet && \bullet \ar[ll] && \bullet \ar[ll] && \bullet \ar[ll]
&& \bullet \ar[ll] \ar[ur]^{\alpha} && \bullet \ar[ll]^{\gamma}
\ar[rr] && \bullet \ar[rr] && \bullet}$$
with the relations $\alpha\beta=0$, $\beta\gamma=0$ and
$\gamma\alpha=0$. Its AR-quiver is given in Fig.~1 below, in which
the projective modules are identified with circles and the
injective modules are identified with squares. The left part $\la$
has two clearly identified connected components (compare with
Lemma \ref{AC03 (1.6)}). Both ends are identified along the
vertical dotted lines, in the inverse order like a Mobiüs band.
Finally, the black diamonds represent the (indecomposable)
Ext-injective modules in $\add\la$.
\begin{figure}[!htb]
\begin{center}
\begin{picture}(0,0)%
\includegraphics{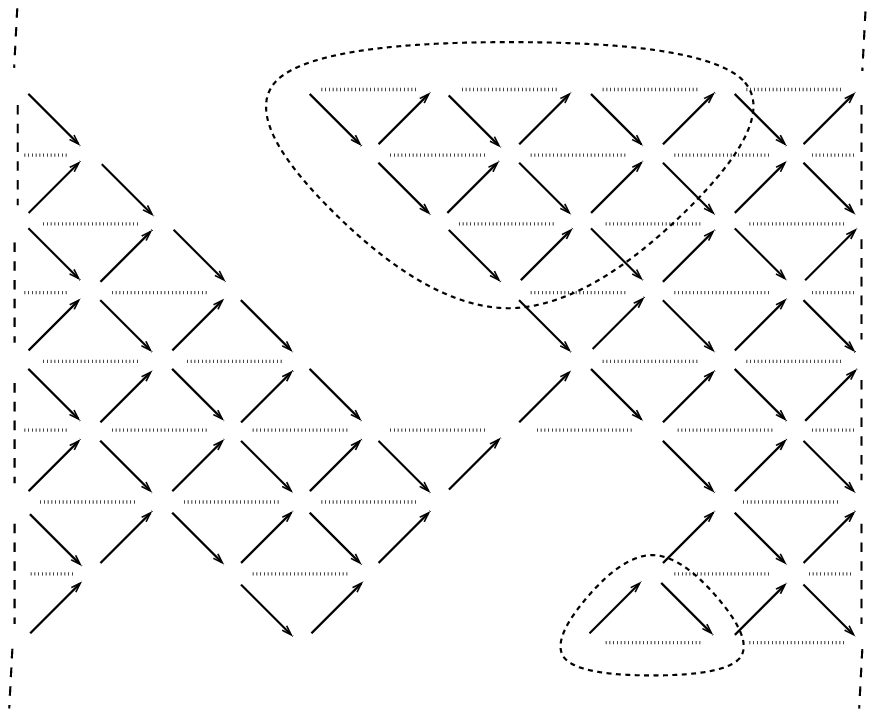}%
\end{picture}%
\setlength{\unitlength}{789sp}%
\begingroup\makeatletter\ifx\SetFigFont\undefined%
\gdef\SetFigFont#1#2#3#4#5{%
  \reset@font\fontsize{#1}{#2pt}%
  \fontfamily{#3}\fontseries{#4}\fontshape{#5}%
  \selectfont}%
\fi\endgroup%
\begin{picture}(20838,17701)(1426,-24920)
\put(8101,-23236){\makebox(0,0)[lb]{\smash{\SetFigFont{6}{7.2}{\rmdefault}{\mddefault}{\updefault}$\xymatrix@1{*[o][F]{\fbox{$\bullet$}}}$}}}
\put(1726,-13111){\makebox(0,0)[lb]{\smash{\SetFigFont{6}{7.2}{\rmdefault}{\mddefault}{\updefault}$\bullet$}}}
\put(3376,-14836){\makebox(0,0)[lb]{\smash{\SetFigFont{6}{7.2}{\rmdefault}{\mddefault}{\updefault}$\bullet$}}}
\put(10126,-18211){\makebox(0,0)[lb]{\smash{\SetFigFont{6}{7.2}{\rmdefault}{\mddefault}{\updefault}$\bullet$}}}
\put(6826,-18211){\makebox(0,0)[lb]{\smash{\SetFigFont{6}{7.2}{\rmdefault}{\mddefault}{\updefault}$\bullet$}}}
\put(8476,-19861){\makebox(0,0)[lb]{\smash{\SetFigFont{6}{7.2}{\rmdefault}{\mddefault}{\updefault}$\bullet$}}}
\put(6826,-21586){\makebox(0,0)[lb]{\smash{\SetFigFont{6}{7.2}{\rmdefault}{\mddefault}{\updefault}$\bullet$}}}
\put(5101,-19861){\makebox(0,0)[lb]{\smash{\SetFigFont{6}{7.2}{\rmdefault}{\mddefault}{\updefault}$\bullet$}}}
\put(3451,-18211){\makebox(0,0)[lb]{\smash{\SetFigFont{6}{7.2}{\rmdefault}{\mddefault}{\updefault}$\bullet$}}}
\put(1726,-16486){\makebox(0,0)[lb]{\smash{\SetFigFont{6}{7.2}{\rmdefault}{\mddefault}{\updefault}$\bullet$}}}
\put(1726,-19936){\makebox(0,0)[lb]{\smash{\SetFigFont{6}{7.2}{\rmdefault}{\mddefault}{\updefault}$\bullet$}}}
\put(3451,-21661){\makebox(0,0)[lb]{\smash{\SetFigFont{6}{7.2}{\rmdefault}{\mddefault}{\updefault}$\bullet$}}}
\put(5101,-16486){\makebox(0,0)[lb]{\smash{\SetFigFont{6}{7.2}{\rmdefault}{\mddefault}{\updefault}$\bullet$}}}
\put(20401,-11536){\makebox(0,0)[lb]{\smash{\SetFigFont{6}{7.2}{\rmdefault}{\mddefault}{\updefault}$\bullet$}}}
\put(21976,-13111){\makebox(0,0)[lb]{\smash{\SetFigFont{6}{7.2}{\rmdefault}{\mddefault}{\updefault}$\bullet$}}}
\put(21976,-16486){\makebox(0,0)[lb]{\smash{\SetFigFont{6}{7.2}{\rmdefault}{\mddefault}{\updefault}$\bullet$}}}
\put(20326,-14836){\makebox(0,0)[lb]{\smash{\SetFigFont{6}{7.2}{\rmdefault}{\mddefault}{\updefault}$\bullet$}}}
\put(18601,-16486){\makebox(0,0)[lb]{\smash{\SetFigFont{6}{7.2}{\rmdefault}{\mddefault}{\updefault}$\bullet$}}}
\put(20326,-18211){\makebox(0,0)[lb]{\smash{\SetFigFont{6}{7.2}{\rmdefault}{\mddefault}{\updefault}$\bullet$}}}
\put(16951,-18211){\makebox(0,0)[lb]{\smash{\SetFigFont{6}{7.2}{\rmdefault}{\mddefault}{\updefault}$\bullet$}}}
\put(21976,-19861){\makebox(0,0)[lb]{\smash{\SetFigFont{6}{7.2}{\rmdefault}{\mddefault}{\updefault}$\bullet$}}}
\put(20326,-21661){\makebox(0,0)[lb]{\smash{\SetFigFont{6}{7.2}{\rmdefault}{\mddefault}{\updefault}$\bullet$}}}
\put(13501,-18136){\makebox(0,0)[lb]{\smash{\SetFigFont{6}{7.2}{\rmdefault}{\mddefault}{\updefault}$\bullet$}}}
\put(3451,-11536){\makebox(0,0)[lb]{\smash{\SetFigFont{6}{7.2}{\rmdefault}{\mddefault}{\updefault}$\fbox{$\bullet$}$}}}
\put(5176,-13111){\makebox(0,0)[lb]{\smash{\SetFigFont{6}{7.2}{\rmdefault}{\mddefault}{\updefault}$\fbox{$\bullet$}$}}}
\put(8401,-16486){\makebox(0,0)[lb]{\smash{\SetFigFont{6}{7.2}{\rmdefault}{\mddefault}{\updefault}$\fbox{$\bullet$}$}}}
\put(6751,-14836){\makebox(0,0)[lb]{\smash{\SetFigFont{6}{7.2}{\rmdefault}{\mddefault}{\updefault}$\fbox{$\bullet$}$}}}
\put(11701,-19861){\makebox(0,0)[lb]{\smash{\SetFigFont{6}{7.2}{\rmdefault}{\mddefault}{\updefault}$\fbox{$\bullet$}$}}}
\put(10051,-21586){\makebox(0,0)[lb]{\smash{\SetFigFont{6}{7.2}{\rmdefault}{\mddefault}{\updefault}$\fbox{$\bullet$}$}}}
\put(18301,-19861){\makebox(0,0)[lb]{\smash{\SetFigFont{6}{7.2}{\rmdefault}{\mddefault}{\updefault}$\xymatrix@1{*+++++[o][F-]{\bullet}}$}}}
\put(16576,-21586){\makebox(0,0)[lb]{\smash{\SetFigFont{6}{7.2}{\rmdefault}{\mddefault}{\updefault}$\xymatrix@1{*+++++[o][F-]{\blacklozenge}}$}}}
\put(21826,-9886){\makebox(0,0)[lb]{\smash{\SetFigFont{6}{7.2}{\rmdefault}{\mddefault}{\updefault}$\fbox{$\bullet$}$}}}
\put(1576,-9886){\makebox(0,0)[lb]{\smash{\SetFigFont{6}{7.2}{\rmdefault}{\mddefault}{\updefault}$\fbox{$\bullet$}$}}}
\put(12976,-7711){\makebox(0,0)[lb]{\smash{\SetFigFont{9}{10.8}{\rmdefault}{\mddefault}{\updefault}$\mathcal{L}_A$}}}
\put(16426,-24736){\makebox(0,0)[lb]{\smash{\SetFigFont{9}{10.8}{\rmdefault}{\mddefault}{\updefault}$\mathcal{L}_A$}}}
\put(18601,-23311){\makebox(0,0)[lb]{\smash{\SetFigFont{6}{7.2}{\rmdefault}{\mddefault}{\updefault}$\blacklozenge$}}}
\put(14851,-23236){\makebox(0,0)[lb]{\smash{\SetFigFont{6}{7.2}{\rmdefault}{\mddefault}{\updefault}$\xymatrix@1{*+++++[o][F-]{\bullet}}$}}}
\put(1426,-23161){\makebox(0,0)[lb]{\smash{\SetFigFont{6}{7.2}{\rmdefault}{\mddefault}{\updefault}$\fbox{$\bullet$}$}}}
\put(21826,-23236){\makebox(0,0)[lb]{\smash{\SetFigFont{6}{7.2}{\rmdefault}{\mddefault}{\updefault}$\fbox{$\bullet$}$}}}
\put(16876,-11536){\makebox(0,0)[lb]{\smash{\SetFigFont{6}{7.2}{\rmdefault}{\mddefault}{\updefault}$\blacklozenge$}}}
\put(15226,-13111){\makebox(0,0)[lb]{\smash{\SetFigFont{6}{7.2}{\rmdefault}{\mddefault}{\updefault}$\blacklozenge$}}}
\put(14926,-16411){\makebox(0,0)[lb]{\smash{\SetFigFont{6}{7.2}{\rmdefault}{\mddefault}{\updefault}$\xymatrix@1{*+++++[o][F-]{\bullet}}$}}}
\put(18676,-13111){\makebox(0,0)[lb]{\smash{\SetFigFont{6}{7.2}{\rmdefault}{\mddefault}{\updefault}$\bullet$}}}
\put(16876,-14761){\makebox(0,0)[lb]{\smash{\SetFigFont{6}{7.2}{\rmdefault}{\mddefault}{\updefault}$\bullet$}}}
\put(13201,-14761){\makebox(0,0)[lb]{\smash{\SetFigFont{6}{7.2}{\rmdefault}{\mddefault}{\updefault}$\xymatrix@1{*+++++[o][F-]{\blacklozenge}}$}}}
\put(15226,-9886){\makebox(0,0)[lb]{\smash{\SetFigFont{6}{7.2}{\rmdefault}{\mddefault}{\updefault}$\bullet$}}}
\put(13501,-11536){\makebox(0,0)[lb]{\smash{\SetFigFont{6}{7.2}{\rmdefault}{\mddefault}{\updefault}$\bullet$}}}
\put(11476,-13111){\makebox(0,0)[lb]{\smash{\SetFigFont{6}{7.2}{\rmdefault}{\mddefault}{\updefault}$\xymatrix@1{*+++++[o][F-]{\bullet}}$}}}
\put(11851,-9886){\makebox(0,0)[lb]{\smash{\SetFigFont{6}{7.2}{\rmdefault}{\mddefault}{\updefault}$\bullet$}}}
\put(9901,-11461){\makebox(0,0)[lb]{\smash{\SetFigFont{6}{7.2}{\rmdefault}{\mddefault}{\updefault}$\xymatrix@1{*+++++[o][F-]{\bullet}}$}}}
\put(8176,-9886){\makebox(0,0)[lb]{\smash{\SetFigFont{6}{7.2}{\rmdefault}{\mddefault}{\updefault}$\xymatrix@1{*+++++[o][F-]{\bullet}}$}}}
\put(18676,-9886){\makebox(0,0)[lb]{\smash{\SetFigFont{6}{7.2}{\rmdefault}{\mddefault}{\updefault}$\blacklozenge$}}}
\end{picture}
    \caption{AR-quiver of the algebra of Example \ref{ex1}}
\end{center}
\end{figure}
\end{my_ex}

Let $A$ be an artin algebra and $\overline{P}$ denote the direct
sum of all indecomposable projective modules in $\la$. In
\cite{ACT04, Sk03}, the algebra $A_\lambda =
\End_A(\overline{P})^{op}$, called the \textit{left support
algebra} of $A$, was studied and shown to be a direct product of
quasitilted algebras. In the above example, one can observe that
$A_\lambda$ is a direct product of (two) hereditary algebras, and
also that $\E_1=\emptyset$ since, equivalently, $\la$ contains no
injective module.  Also, the left part is given by the modules
which are not successors of any injective module. This is not a
coincidence as the following results show.

\begin{proposition}\label{prop hered}%
Let $A$ be an algebra of Gorenstein dimension at most one.  The
left support algebra $A_\lambda$ is a direct product of hereditary
algebras.
\end{proposition}
\begin{proof}
Since $\la\subseteq \mod A_{\lambda}$ by \cite{ACT04}, it suffices
to show that if $P$ is a projective module in $\la$ and
$\xymatrix@1@C=15pt{M\ar[r] & P}$ is an irreducible morphism, then
$M$ is projective. If $M$ is not projective, then $\tau M\neq 0$
and thus $\Hom_A(\tau^{-1}(\tau M), P)\neq 0$.  By Lemma \ref{lem
pd1}, this gives $\id_A\tau M>1$.  The Gorenstein property then
implies $\pd_A\tau M>1$, a contradiction to $\tau M\in \add\la$.
Thus $M$ is projective.
\end{proof}

\begin{my_cor}\label{cor E_1}%
Let $\C_H$ be a cluster category, $T$ be a tilting object in
$\C_H$ and $A=\End_{\C_H}(T)^{op}$. If $A$ is not hereditary, then
$\E_1=\emptyset$.
\end{my_cor}
\begin{proof}
We first show that if $I'$ is an injective module in $\la$ and
$\mor{f}{I'}{M}$ is an irreducible morphism, with $M$
indecomposable, then $M\in \la$.  Indeed, assume that
$M\notin\la$. Then $\tau M$ is Ext-injective in $\la$ and $\tau
M\in\E_1\cup\E_2$ (observe that $M$ is not projective since
$\Ext^1_{C_H}(T,T)=0$). If $M\in\E_1$, then there exists an
injective module $I''$ in $\la$ together with a path
$\xymatrix@1@C=15pt{I''=N_0 \ar[r] & \cdots \ar[r] &N_m=\tau M
\ar[r] &I'}$
in $\la$. Now, since $A_\lambda$ is hereditary by Proposition
\ref{prop hered}, $\tau M$ is injective, and so $M=0$, a
contradiction.  Hence $\tau M\in\E_2$. Then, there exists an
indecomposable projective module $P\notin\la$ and a sectional path
$\path{\delta}{P}{M}$.  Let $T'$ be the direct summand of $T$
corresponding to $P$ and $T''$ be the direct summand of $T$ such
that $T''[2]$ corresponds to $I'$. Then, lifting the path $\delta$
in $\C_H$, and using the fact that this path does not factor
through $I'$ (since $I'\in\la$), yields a sectional path from $T'$
to $T''[1]$, a contradiction to $\Hom_{\C_H}(T, T[1])\neq 0$.
Therefore $M\in\la$.

Now, assume that $I$ is an injective module in $\la$.  Let $\Ga$
be the connected component of the AR-quiver of $A$ containing $I$
and $\Sigma$ be the maximal full, connected and convex subquiver
of $\Ga$ containing only indecomposable injective modules,
including $I$. Observe that since $\la$ is closed under
predecessors, and in view of the first part of the proof, any
injective module in $\Sigma$ lies in $\la$.  Now, dualizing the
arguments in the proof of Lemma \ref{lem nsr} yields an injective
module $I'$ in $\Sigma$ together with an irreducible morphism
$\xymatrix@1@C=15pt{I'\ar[r]&M}$, where $M$ is not injective. But
since $I'\in\la$, we get $M\in \la$ by the first part of the
proof, a contradiction to the fact that $A_\lambda$ is a direct
product of hereditary algebras (since $M$ is not injective).
\end{proof}

Thus, the left part of a cluster-tilted not hereditary algebra
contains no injective module. We get the following easy
consequence of Lemma \ref{AC03 (1.6)}.

\begin{my_cor}\label{cor la}%
Let $A$ be a cluster-tilted algebra. If $A$
is not hereditary, then
$$\begin{array}{c} \la=\{M\in\ind A : M \mbox{ is not a successor
of an injective module}\}\\
\ra=\{M\in\ind A : M \mbox{ is not a predecessor of a
projective module}\} %
\end{array}$$
\end{my_cor}

The following lemma, whose proof follows directly from the above
corollary, will be useful in the next section.

\begin{my_lem}\label{lem LT}%
Let $\C_H$ be a cluster category, $T$ be a
tilting object in $\C_H$ and $A=\End_{\C_H}(T)^{op}$. The functor
$\Hom_{\C_H}(T,-)$ induces an equivalence
$\xymatrix@1@C=15pt{\mathcal{L}_T \ar[r] & \la},$
where $\mathcal{L}_T$ is the set of all indecomposable objects $M$
in $\C_H\setminus \add T[1]$ such that if there exists a path from
an indecomposable object in $\add T[2]$ to $M$, then (at least one
morphism in) this path factors through $\add T[1]$.
\end{my_lem}

Proposition \ref{prop hered} has another nice direct consequence.

\begin{my_cor}\label{cor hered}%
Let $A$ be an algebra of Gorenstein dimension at most one. Then
$A$ is hereditary if and only if $A\in\add\la$, and this occurs if
and only if $A$ is quasitilted (see \cite[(II.1.14)]{HRS96}).
\end{my_cor}

%
%
\subsection{Endomorphism algebras of $\la$-slices}
    \label{la-slices}

Here, we introduce the concept of $\la$-slices and show that if
$A$ is cluster-tilted, then these $\la$-slices induce tilting
modules whose endomorphism algebras are again cluster-tilted.

We first recall the following definition: let $(\Ga, \tau)$ be a
connected translation quiver. A connected full subquiver $\Sigma$
of $\Ga$ is a \textit{section} in $\Ga$ if:
\begin{enumerate}
\item[(S1)] $\Sigma$ is acyclic;
\item[(S2)] For each $x\in\Sigma$, there exists a unique $n\in\mathbb{Z}$ such that $\tau^nx\in\Sigma$;
\item[(S3)] $\Sigma$ is convex in $\Ga$.
\end{enumerate}

This definition was motivated by the study of tilted algebras. The
well-known criterion of Liu and Skowro{\'n}ski (see
\cite[(VIII.5.6)]{ASS06} for instance)
asserts that an algebra $A$ is tilted if and only if its
AR-quiver has a connected component containing a faithful section
$\Sigma$ such that $\Hom_A(X, \tau Y)=0$ for each $X,Y\in \Sigma$.
These faithful sections were called \textit{complete slices in
$\mod A$}.

By \cite[(Theorem B)]{ACT04}, an algebra $A$ is left supported if
and only if each connected component of $A_\lambda$ is a tilted
algebra and the restriction of $E$ (see Section \ref{Left and
right supported algebras}) to this component is a complete slice.
Since, by construction, we have $\la\subseteq \mod
A_\lambda\subseteq \mod A$, this motivates the following
definition:

\begin{my_def}\label{defn la-slice}%
Let $A$ be an algebra and $A_\lambda=A_1\times\cdots \times A_m$
be its left support algebra. An \emph{$\la$-slice} is a direct
product $S=S_1 \times \cdots \times S_m$, with each $S_i$ a
complete slice in $\mod A_i\cap \la$.
\end{my_def}

Such $\la$-slices do not always exist, for instance when
$A=A_\lambda$ is a quasitilted not tilted algebra, or worse when
$\la=\emptyset$. Here, we give two canonical examples of
$\la$-slices when $A$ is cluster-tilted.

\begin{my_ex}\label{ex2}%
Let $A$ be a cluster-tilted algebra such that $\la\neq\emptyset$.
\begin{enumerate}
\item[(a)] By Proposition \ref{prop hered}, $A_\lambda$ is a direct product of hereditary algebras. Then, the
full subquiver generated by the set $\Sigma_P=\{P_1, \dots, P_n\}$
of indecomposable projective modules in $\la$ is an $\la$-slice.
\item[(b)] By Proposition \ref{prop ls}, $A$ is left supported.  Hence, by the above mentioned theorem, the
direct sum $E$ of the indecomposable Ext-injective modules in
$\add\la$ is an $\la$-slice (compare with Example \ref{ex1}).
\end{enumerate}
\end{my_ex}

Clearly, these two examples are extremal, in the sense that any
$\la$-slice lies between these two.  Morever, we get the
following:

\begin{my_lem}\label{lem exch}%
Let $A$ be a cluster-tilted not hereditary algebra.  Let
$\Sigma_P$ be the $\la$-slice generated by the projective modules
in $\la$. Then any $\la$-slice $\Sigma$ can be reached from
$\Sigma_P$ by a finite number of almost split exchanges.
\end{my_lem}
\begin{proof}
Let $\Sigma$ be an $\la$-slice and $P_1, \dots, P_n$ be the
vertices of $\Sigma_P$. Assume that $\Sigma_P$ has a source $P_i$
which is not in $\Sigma$. Then, replacing in $\Sigma_P$ the module
$P_i$ by $\tau^{-1}P_i$ and all arrows $\xymatrix@1@C=15pt{P_i
\ar[r] & P_j}$ by their corresponding arrows
$\xymatrix@1@C=15pt{P_j \ar[r] & \tau^{-1}P_i}$ yields a new
$\la$-slice $\Sigma'_P$.  By iterating this procedure and invoking
that $\la$ is finite by Proposition \ref{prop ls}, we get after
finitely many steps the $\la$-slice $\Sigma$.
\end{proof}

Clearly, by using the above procedure, the number of needed almost
split exchanges to reach the $\la$-slice $\Sigma$ is uniquely
determined. Indeed, if $\Sigma=\{S_1, \dots, S_n\}$ with
$S_i=\tau^{-t_i}P_i$ for each $i$, then the number of required
exchanges is given by $t_{\Sigma}=\sum_{i=1}^nt_i$. In particular,
when $\Sigma=E$ (see Section \ref{Left and right supported
algebras}), then $t_{\Sigma}=|\la|-n$, where $n$ denotes the
number of indecomposable projective modules in $\la$.

We can now prove Theorem \ref{Theorem 3}. Observe that, here
again, we keep the same notation for an $A$-module and its
preimage in $\C_H$.

\begin{my_thm2}\label{Theorem3}%
Let $\C_H$ be a cluster category, $T$ be
a tilting object in $\C_H$ and $A$ be the cluster-tilted algebra
$\End_{\C_H}(T)^{op}$. Assume that $A$ is not hereditary and let
$\Sigma$ be an $\la$-slice. Also, let $F=\oplus_{i=1}^m P_i$
denote the direct sum of all indecomposable projective modules not
in $\la$. Then,
\begin{enumerate}
\item[(a)] $T_{\Sigma}=\Sigma\oplus F$ is a tilting
$A$-module;
\item[(b)] The algebra $A_\Sigma=\End_A(T_{\Sigma})^{op}$ is
isomorphic to $\End_{\C_H}(T_{\Sigma})^{op}$.  In particular,
$A_\Sigma$ is cluster-tilted;
\item[(c)] The quiver of $A_\Sigma$ is obtained from that of $A$ with $t_{\Sigma}$ reflections at sinks.
\end{enumerate}
\end{my_thm2}
\begin{proof}
(a). We prove a more general fact.  Let $n$ be the number of
indecomposable projective modules in $\la$ and $S_1, \dots, S_n$
be $A$-modules in $\la (\subseteq \mod A_\lambda)$ such that
$\Hom_{A_{\lambda}}(S_i, \tau S_j)=0$ for all $i,j$.  Since $\la$
is closed under predecessors, we get $0=\Hom_A(S_i, \tau
S_j)=\Ext_A^1(S_j,S_i)$ for all $i,j$. Let
$\Sigma=\oplus_{i=1}^nS_i$ and $T_{\Sigma}=S\oplus F$. Then,
$\Ext_A^1(\Sigma,F)\cong D\Hom_A(F,\tau \Sigma)=0$ an since
$\pd_AT_{\Sigma}\leq 1$ by construction, $T_{\Sigma}$ is a tilting
$A$-module.\\
(b). By Theorem \ref{Theorem 1}, $T_{\Sigma}$ is a tilting object
in $\C_H$. So $\End_{\C_H}(T_{\Sigma})^{op}$ is cluster-tilted. In
view of the equivalence $\C_H/\add T[1]\cong \mod A$, it then
suffices to show that no morphism between two direct summands of
$T_{\Sigma}$ in $\C_H$ factors through $\add T[1]$. We prove it by
induction on number $t_{\Sigma}$ of necessary almost split
exchanges to reach $\Sigma$ from the $\la$-slice $\Sigma_P$
generated by the set of indecomposable projective modules in $\la$
(see Lemma \ref{lem exch}). Let $\Sigma=\{S_1, \dots, S_n\}$ and
$T_1,\dots, T_m$ be the indecomposable direct summands of $T$
corresponding to the indecomposable direct summands $P_1, \dots,
P_m$ of $F$. If $t_{\Sigma}=0$, then $\Sigma=\Sigma_P$ and
$T_{\Sigma}=T$. The claim then follows from $\Hom_{\C_H}(T,
T[1])=0$. Assume that $t_{\Sigma}>0$. Since each connected
component of $\Sigma$ is acyclic, $\Sigma$ contains some sinks.
Also, since $A_\lambda$ is hereditary, some of these sinks are not
projective. Assume that $S_1$ is a non-projective sink in $\Sigma$
and consider the $\la$-slice $\Sigma'$ obtained by replacing in
$\Sigma$ the module $S_1$ by $\tau S_1$ and all arrows
$\xymatrix@1@C=15pt{S_i \ar[r] & S_1}$ by their corresponding
arrows $\xymatrix@1@C=15pt{\tau S_1 \ar[r] & S_i}$. So
$\Sigma'=\{\tau S_1, S_2, \dots, S_n\}$. We have
$t_{\Sigma'}<t_{\Sigma}$, and thus by induction no morphism
between two direct summands of $T_{\Sigma'}=\Sigma'\oplus F$ in
$\C_H$ factors through $\add T[1]$.

To prove our claim, we then have to show that no morphism in one
of the Hom-spaces: (i) $\Hom_{\C_H}(S_1, S_i)$, (ii)
$\Hom_{\C_H}(S_i, S_1)$, (iii) $\Hom_{\C_H}(S_1, T_j)$ and (iv)
$\Hom_{\C_H}(T_j, S_1)$, for $2\leq i\leq n$ and $1\leq j\leq m$,
factors through $\add T[1]$.
\begin{enumerate}
\item[(i)] For each $i=2, \dots, n$, we have
$\Hom_{\C_H}(S_1, S_i)\cong \Hom_{\C_H}(\tau S_1, \tau S_i)=0$
because $T_{\Sigma'}$ is a tilting object in $\C_H$. This is
sufficient.
\item[(ii)] Let $\xymatrix@1@C=15pt{0\ar[r] & \tau S_1 \ar[r] &
\oplus_{k=1}^{q}S_{1,k} \ar[r] & S_1 \ar[r] & 0}$ be the short
exact sequence ending in $S_1$.  Since $S_1$ is a sink in
$\Sigma$, it follows from the construction of $\Sigma$ by
$\Sigma_P$ (see Lemma \ref{lem exch}) that $S_{1,k}$ is a vertex
in $\Sigma\cap\Sigma'$ for each $k$.  Now, if $\mor{f}{S_i}{S_1}$
is a morphism in $\C_H$, then $f$ lifts to a morphism from $S_i$
to $\oplus_{k=1}^{q}S_{1,k}$.  By the induction hypothesis, this
morphism does not factor through $\add T[1]$, and hence the same
holds for $f$.
\item[(iii)] As in (i), for $j=1, \dots, m$, we have
$\Hom_{\C_H}(S_1, T_j)\cong \Hom_{\C_H}(\tau S_1, \tau T_j)=0$
because $T_{\Sigma'}$ is a tilting object in $\C_H$. This is
sufficient.
\item[(iv)] Finally, since $\Hom_{\C_H}(T,T[1])=0$, no morphism from some $T_j$ to $S_1$ factors
through $\add T[1]$.
\end{enumerate}
Consequently,
$\End_A(T_{\Sigma})^{op}\cong\End_{\C_H}(T_{\Sigma})^{op}$ is
cluster-tilted.\\
(c). We first recall a general fact : let $A=\End_{\C_H}(T)^{op}$
be a cluster-tilted algebra. Also, let $T=\overline{T}\oplus M$,
with $M$ indecomposable, and $M^*$ be the other complement for
$\overline{T}$. Finally, let $T^*=\overline{T}\oplus M^*$ and
$A^*=\End_{\C_H}(T^*)^{op}$. By result of Buan, Marsh and Reiten
\cite{BMR05} the quivers $Q_A$ of $A$ and $Q_{A^*}$ of $A^*$ are
related by the quiver mutation formula of Fomin and Zelevinsky. In
particular, when $M$ corresponds to a sink in $Q_A$, then
$Q_{A^*}$ is obtained from $Q_A$ by performing a reflection at
this sink.

In our case, because $A_{\lambda}$ is hereditary, each almost
split exchange performed in the proof of Lemma \ref{lem exch} (in
order to reach $\Sigma$ from $\Sigma_P$) coincides in $\C_H$ with
an almost split exchange of an indecomposable direct summand $M$
of a certain tilting object, say
$T_{\Sigma'}=\overline{T_{\Sigma'}}\oplus M$, by the other
complement $M^*=\tau^{-1}M$ of $\overline{T_{\Sigma'}}$ (see
Remark \ref{rem ase}). Moreover, $M$ corresponds to a sink in the
quiver associated with
$\End_A(T_{\Sigma'})^{op}\cong\End_{\C_H}(T_{\Sigma'})^{op}$.
Therefore, by \cite{BMR05}, this almost split exchange coincides
with a reflection at a sink in the quiver of
$\End_A(T_{\Sigma'})^{op}$. Now, since, in the notations of (b),
$A_{\Sigma_P}=A$ and $\Sigma$ can be reached from $\Sigma_P$ with
$t_{\Sigma}$ almost split exchanges, this means that the quiver of
$A_{\Sigma}$ can be obtained from that of $A$ by performing
$t_{\Sigma}$ reflections at sinks.
\end{proof}

Recall from Theorem \ref{thm ACT} that $A$ is left supported if
and only if the $A$-modules $L=E\oplus F$ and $U=E_1\oplus
E_2\oplus F$ are tilting modules. Since $L$ is induced by the
Ext-injective modules in $\add\la$, it follows from the above
theorem that $\End_A(L)^{op}$ is cluster-tilted.  We now show that
the same holds for $\End_A(U)^{op}$ although $U$ does not arise
from an $\la$-slice. At this point, we stress that since $E_1=0$
by Corollary \ref{cor E_1}, we have $U=\tau^{-1}E_2\oplus
F=\tau^{-1}E\oplus F$.

We need the following lemma (compare with Example \ref{ex1}).

\begin{my_lem}\label{lem source}%
Let $A$ be an algebra and $\E$ be the set of all indecomposable
Ext-injective modules in $\add\la$. If $M$ is a source in $\E$ and
$\mor{f}{M}{N}$ is an irreducible morphism, with $N$
indecomposable, then $N\in\E$ or $N$ is projective.
\end{my_lem}
\begin{proof}
Indeed, if $N\notin \E$ and $N$ is not projective, then $\tau N$
exists and belongs to $\la$ (since it is a predecessor of $M$).
Moreover, $N\notin \E$ implies $N\notin\la$ since $\E$ is closed
under successors in $\la$ by \cite[(3.4)]{ACT04}. So $\tau N\in
\E$. But this contradicts the fact that $M$ is a source in $\E$.
So $N\in\E$ or $N$ is projective.
\end{proof}

\begin{proposition}\label{prop U}%
Let $A$ be a cluster-tilted algebra
which is not hereditary and $U=\tau^{-1}E\oplus F$ be as above.
Then,
\begin{enumerate}
\item[(a)] $U$ is a tilting
$A$-module;
\item[(b)] The algebra $A_U=\End_A(U)^{op}$ is
isomorphic to $\End_{\C_H}(U)^{op}$.  In particular, $A_U$ is
cluster-tilted;
\item[(c)] The quiver of $A_U$ is obtained from that of $A$ with $|\la|$ reflections at sinks.
\end{enumerate}
\end{proposition}
\begin{proof}
(a). This follows from Theorem \ref{thm ACT}.\\ %
(b) and (c). By Theorem \ref{Theorem 1}, $U$ is a tilting object
in $\C_H$. Also, by continuing the procedure in the proof of Lemma
\ref{lem exch}, $\tau^{-1}E$ is obtained from $E$ by performing
$n$ almost split exchanges in $\mod A$, where $n$ denotes the
number of projective modules in $\la$. By Lemma \ref{lem source}
and Remark \ref{rem ase}, these exchanges correspond in $\C_H$ to
(almost split) exchanges of tilting object. So, the quiver of
$\End_{\C_H}(U)^{op}$ is obtained from that of
$\End_{\C_H}(L)^{op}$ with $n$ reflections at sinks. Also, as in
the proof of Theorem \ref{Theorem3}, one can show by induction
that $\End_A(U)^{op}\cong\End_{\C_H}(U)^{op}$. Since, by Theorem
\ref{Theorem3}, the quiver of $\End_{\C_H}(L)^{op}$ is obtained
from that of $A$ with $|\la|-n$ reflections at sinks, this proves
(c).
\end{proof}

\begin{my_ex}\label{ex3}%
Let $A$ be the cluster-tilted not hereditary algebra of Example
\ref{ex1}.  Let $E$ be the direct sum of the indecomposable
Ext-injective modules in $\add\la$ (these identified with black
diamonds) and $F$ be the direct sum of the three indecomposable
projective modules not lying in $\la$.  As usual, let $L=E\oplus
F$ and $U=\tau^{-1}E\oplus F$.
\begin{enumerate}
\item[(a)] The algebra $\End_A(L)^{op}$ is the cluster-tilted algebra given by the quiver
$$\xymatrix@1@C=8pt@R=8pt{&&&&&&&&& \bullet \ar[dr]^{\beta}\\ %
\bullet \ar[rr] && \bullet \ar[rr] && \bullet \ar[rr] && \bullet
&& \bullet \ar[ll] \ar[ur]^{\alpha} && \bullet \ar[ll]^{\gamma}
\ar[rr] && \bullet && \bullet \ar[ll]}$$
with the relations $\alpha\beta=0$, $\beta\gamma=0$ and
$\gamma\alpha=0$.
\item[(b)] The algebra $\End_A(U)^{op}$ is the cluster-tilted algebra given by the quiver
$$\xymatrix@1@C=8pt@R=8pt{&&&&&&&&& \bullet \ar[dr]^{\beta}\\ %
\bullet \ar[rr] && \bullet \ar[rr] && \bullet \ar[rr] && \bullet
\ar[rr] && \bullet \ar[ur]^{\alpha} && \bullet \ar[ll]^{\gamma} &&
\bullet \ar[ll] && \bullet \ar[ll]}$$
with the relations $\alpha\beta=0$, $\beta\gamma=0$ and
$\gamma\alpha=0$.
\end{enumerate}
\end{my_ex}

In the above example, one can observe that the quiver of the
algebra $A_U=\End_A(U)^{op}$ has no sink, meaning that
$\mathcal{L}_{A_U}=\emptyset$. The following two results explain
this phenomenon.  Here, the notation $\mathcal{L}_T$ refers to the
subcategory of $\C_H$ introduced in Lemma \ref{lem LT} and
$\mathcal{R}_T$ refers to its analogue for the right part.

\begin{proposition}\label{prop lara}%
Let $\C_H$ be a cluster category, $T$ be a tilting object in
$\C_H$ and $A=\End_{\C_H}(T)^{op}$ be cluster-tilted not
hereditary.  Assume that $\Sigma=\{S_1,\dots, S_n\}$ is an
$\la$-slice having a source $S_1$ such that $\tau^{-1}S_1\in\la$.
Let $\Sigma'=\{\tau^{-1}S_1, S_2, \dots, S_n\}$ be the $\la$-slice
obtained from $\Sigma$ by performing an almost split exchange at
$S_1$. Let $T_\Sigma=\Sigma\oplus F$ and
$T_{\Sigma'}=\Sigma'\oplus F$. Then, in $\C_H$,
\begin{enumerate}
\item[(a)] $\LL_{T_{\Sigma'}}= \LL_{T_{\Sigma}}\setminus \{S_1\}$.
\item[(b)] $\RR_{T_{\Sigma'}}=\RR_{T_{\Sigma}} \cup \{\tau S_1\}$.
\end{enumerate}
In particular, $|\LL_{T_{\Sigma'}}|+|\RR_{T_{\Sigma'}}|=
|\LL_{T_{\Sigma}}|+|\RR_{T_{\Sigma}}|$.
\end{proposition}
\begin{proof}
We only prove (a) since the proof of (b) is dual. \\
(a). Let $S'_1=\tau^{-1}S_1$ and $\overline{T}=
T_\Sigma\setminus\{S_1\}=T_{\Sigma'}\setminus\{S'_1\}$.
We first prove that
$\LL_{T_{\Sigma'}}\subseteq\LL_{T_{\Sigma}}\setminus\{S_1\}$.  Let
$M\in\LL_{T_{\Sigma'}}$ and assume that $M\notin\LL_{T_{\Sigma}}$.
Hence, there exists an indecomposable direct summand $T'$ of
$T_{\Sigma}$ together with a path in $\C_H$ of the form
$\path{\delta}{T'[2]}{M}$ which does not factor through $\add
T_{\Sigma}[1]$.
\begin{enumerate}
\item[(i)] If $T'\in\add\overline{T}$, then $T'\in\add T_{\Sigma'}$
and it follows from $M\in\LL_{T_{\Sigma'}}$ that $\delta$ factors
through $\add S'_1[1]=\add S_1$. But this gives a path from $T'[2]$
to $S_1$, contradicting the fact that $S_1\in\LL_{T_{\Sigma}}$ by
Lemma \ref{lem LT}.
\item[(ii)] Let $T'=S_1$.  Since $S_1$ is a source in
$\Sigma$, the path $\path{\delta}{S_1[2]}{M}$ factors through $\add
\overline{T}[2]$, are we are back to the situation (i).
\end{enumerate}
Therefore, $\LL_{T_{\Sigma'}}\subseteq \LL_{T_{\Sigma}}$, and
since $S_1=S'_1[1]\notin\LL_{T_{\Sigma'}}$, we get
$\LL_{T_{\Sigma'}}\subseteq\LL_{T_{\Sigma}}\setminus\{S_1\}$.
We now prove the inverse inclusion. Let
$M\in\LL_{T_{\Sigma}}\setminus\{S_1\}$ and assume that there is a
path $\path{\delta}{T'[2]}{M}$ in $\C_H$, for some indecomposable
direct summand $T'$ of $T_{\Sigma'}$.  We need to show that
$\delta$ factors through $\add T_{\Sigma'}[1]$.
\begin{enumerate}
\item[(i)] Assume that $T'=S_1'$.  Since $S'_1[2]=S_1[1]$ and
$S_1$ is a source in $\Sigma$, the path $\delta$ factors through
$\add\overline{T}[1]$ and so factors through $\add
T_{\Sigma'}[1]$.
\item[(ii)] If $T'\neq S_1'$, then $T'\in\add\overline{T}$,
and since $M\in\LL_{T_{\Sigma}}$, the path $\delta$ factors
through $\add T_{\Sigma}[1]$. But then $\delta$ factors through
$\add T_{\Sigma'}[1]$ by (i).
\end{enumerate}
\end{proof}

\begin{my_cor}\label{cor lara}%
Let $\C_H$ be a cluster category, $T$ be a
tilting object in $\C_H$ and $A=\End_{\C_H}(T)^{op}$. Assume that
$A$ is not hereditary and let $\Sigma_P=\{P_1, \dots, P_n\}$ be
the $\la$-slice generated by the indecomposable projective modules
in $\la$. Also, let $\Sigma=\{S_1, \dots, S_n\}$ be $\tau^{-1}E$
or an $\la$-slice, and assume that $\Sigma$ can be reached from
$\Sigma_P$ with $t_\Sigma$ almost split exchanges (as in Lemma
\ref{lem exch}). Finally, let $T_{\Sigma}=\Sigma\oplus F$.
\begin{enumerate}
\item[(a)] $|\LL_{T_{\Sigma}}|=|\LL_T|-t_{\Sigma}$.
\item[(b)] $|\RR_{T_{\Sigma}}|=|\RR_T|+t_{\Sigma}$.
\end{enumerate}
In particular, for $U=\tau^{-1}E\oplus F$, we get $|\LL_{U}|=0$
and $|\RR_{U}|=|\RR_T|+|\LL_T|$.
\end{my_cor}

%
%
\newcommand{\etalchar}[1]{$^{#1}$}
\providecommand{\bysame}{\leavevmode\hbox
to3em{\hrulefill}\thinspace}
\providecommand{\MR}{\relax\ifhmode\unskip\space\fi MR }
\providecommand{\MRhref}[2]{%
  \href{http://www.ams.org/mathscinet-getitem?mr=#1}{#2}
} \providecommand{\href}[2]{#2}

\end{document}